\RequirePackage{fix-cm}
\documentclass[smallextended]{svjour3}       
\smartqed  
\usepackage{graphicx}
\usepackage{todonotes}  
\usepackage{mathtools}
\usepackage{amsmath}
\usepackage{amsfonts}
\usepackage{enumitem}
\usepackage{algorithm}
\usepackage{algorithmic}
\usepackage[algo2e,ruled,vlined]{algorithm2e}
\usepackage{booktabs}  
\usepackage{graphicx}
\usepackage{multirow}
\usepackage{cite}
\usepackage{url}
\def\xx{\boldsymbol{x}}
\def\pp{\boldsymbol{p}}
\def\ss{\boldsymbol{s}}

\newcommand{\B}{\textbf}
\DeclareMathOperator*{\argmin}{argmin}
\DeclareMathOperator{\ReLU}{ReLU}
\DeclareMathOperator{\Sigmoid}{Sigmoid}
\begin{document}

\title{Revisiting local branching with a machine learning lens \thanks{This work is an extension of \cite{Liu_Fischetti_Lodi_2022}}
}


\author{Defeng Liu \and Matteo Fischetti \and Andrea Lodi}

\institute{Defeng Liu \at
             CERC, Polytechnique Montr\'eal, Canada,
              \email{defeng.liu@polymtl.ca }
\and Matteo Fischetti \at  Department of Information Engineering, University of Padua, 
Italy, \email{matteo.fischetti@unipd.it}
\and Andrea Lodi \at
            CERC, Polytechnique Montr\'eal, Canada, and Jacobs Technion-Cornell Institute, Cornell Tech and Technion - IIT, USA,
            \email{andrea.lodi@cornell.edu}
}

\date{Received: date / Accepted: date}

\maketitle

\begin{abstract}
Finding high-quality solutions to mixed-integer linear programming problems (MILPs) is of great importance for many practical applications. In this respect, the refinement heuristic \emph{local branching} (LB) has been proposed to produce improving solutions and has been highly influential for the development of local search methods in MILP. The algorithm iteratively explores a sequence of solution neighborhoods defined by the so-called \emph{local branching constraint}, namely, a linear inequality limiting the distance from a reference solution. For a LB algorithm, the choice of the neighborhood size is critical to performance. 
In this work, we study the relation between the size of the search neighborhood and the behavior of the underlying LB algorithm, and we devise a leaning based framework for predicting the best size for the specific instance to be solved. Furthermore, we have also investigated the relation between the time limit for exploring the LB neighborhood and the actual performance of LB scheme, and devised a strategy for adapting the time limit.
We computationally show that the neighborhood size and time limit can indeed be learned, leading to improved performances and that the overall algorithm generalizes well both with respect to the instance size and, remarkably, across instances.
\end{abstract}

\section{Introduction}
\label{sec:intro}



\emph{Mixed-integer linear programming} (MILP) is a main paradigm for modeling complex combinatorial problems. The exact solution of a MILP model is generally attempted by a branch-and-bound (or branch-and-cut)\cite{land2010automatic} framework. Although state-of-the-art MILP solvers experienced a dramatic performance improvement over the past decades, due to the NP-hardness nature of the problem, the computation load of finding a provable optimal solution for the resulting models can be heavy. In many practical cases, feasible solutions are often required within a very restricted time frame. Hence, one is interested in finding solutions of good quality at the early stage of the computation. In fact, it is also appealing to discover early incumbent solutions in the exact enumerate scheme, which improves the primal bound and reduces the size of the branch-and-bound tree by pruning more nodes \cite{berthold2013measuring}.   

In this respect, the concept of heuristic is well rooted as a principle underlying the search of high-quality solutions. In the literature, a variety of heuristic methods have proven to be remarkably effective, e.g., \emph{local branching} \cite{fischetti2003local}, \emph{feasibility pump} \cite{fischetti2005feasibility}, \emph{RINS} \cite{danna2005exploring}, \emph{RENS} \cite{berthold2014rens}, \emph{proximity search} \cite{fischetti2014proximity}, \emph{large neighborhood search} \cite{gendreau2010handbook}, etc. More details of these developments are reviewed in \cite{berthold2006primal, fischetti2010heuristics, gendreau2010handbook}. In this paper, we focus on local branching, a \emph{refinement heuristic} that iteratively produces improved solutions by exploring suitably predefined solution neighborhoods.

\begin{figure}[htbp!]
    \centering
    \includegraphics[width=0.67\textwidth]{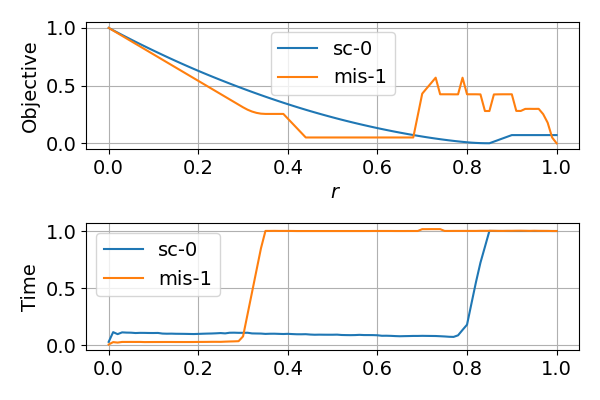}
    \vspace{-4mm}
    \caption{Evaluation of the size of LB neighborhood on a set covering instance (sc-0) and a maximum independent set instance (mis-1). The neighborhood size $k$ is computed as $k = r \times N$, where $N$ is the number of binary variables, and $r \in [0,1]$. A time limit is imposed for each neighborhood exploration. }
    \vspace{-6mm}
    \label{fig:two_examples}
\end{figure}

Local branching (LB) was one of the first methods using a generic MILP solver as a black-box tool inside a heuristic framework. 
Given an initial feasible solution, the method first defines a solution neighborhood through the so-called \emph{local branching constraint}, then explores the resulting subproblem by calling a black-box MILP solver. For a LB algorithm, the choice of neighborhood size is crucial to performance. In the original LB algorithm \cite{fischetti2003local}, the size of neighborhood is mostly initialized by a small value, then adjusted in the subsequent iterations. Although these conservative settings have the advantage of yielding a series of easy-to-solve subproblems, each leading to a small progress of the objective, there is still a lot of space for improvement. As discussed in \cite{fischetti2014proximity}, a significantly better performance can be potentially achieved with an ad-hoc tuning of the size of the neighborhood. Our observation also shows that the ``best" size is strongly dependent on the particular MILP instance. 
To illustrate this, the performance of different LB neighborhood size settings for two MILP instances are compared in Figure \ref{fig:two_examples}. 
In principle, it is desirable to have neighborhoods to be relatively small to allow for an efficient exploration, but still large enough to be effective for finding improved solutions. Nonetheless, it is reasonable to believe that the size of an ideal neighborhood is correlated with the characteristics of the particular problem instance. 

Furthermore, it is worth noting that, in many applications, instances of the same problem are solved repeatedly. Problems of real-world applications have a rich structure. While more and more datasets are collected, patterns and regularities appear. Therefore, problem-specific and task-specific knowledge can be learned from data and applied to the corresponding optimization scenario. This motives a broader paradigm of learning to guide the neighborhood search in refinement heuristics.

In this paper, we investigate a learning framework for sizing the search neighborhood of local branching. In particular, given a MILP instance, we exploit patterns in both the structure of the problem and the information collected from the solving process to predict the size of the LB neighborhood and the time limit for exploring the neighborhood, with the aim of maximizing the performance of the underlying LB algorithm. We computationally show that the neighborhood size and the time limit can indeed be learned, leading to improved performances, and that the overall algorithm generalizes well both with respect to the instance size and, more surprisingly, across instances.

We note that a shorter conference version of this paper appeared in \cite{Liu_Fischetti_Lodi_2022}. Our initial conference paper did not consider the effect of the time limit for solving the LB neighborhood on the performance of overall algorithm. This extended paper addresses this issue and provides additional analysis on applying our refined LB algorithms as a primal heuristic within the MILP sovler.

The paper is organized as follows. In Section \ref{sec:related_work}, we give a review of the related works in the literature. In Section \ref{sec:background}, we introduce the basic local branching scheme and some relevant concepts. In Section \ref{sec:method}, we present our methodology for learning to search in the local branching scheme. In Section \ref{sec:result} describes the setup of our experiments and reports the results to validate our approach. In Section \ref{sec:implementations}, we apply local branching as a primal heuristic and provide two possible implementations of how our local branching scheme interacts with the MILP solver. Section \ref{sec:discussion} concludes the paper and discusses future research.


\section{Related work}
\label{sec:related_work}
Recently, the progress in machine learning (ML) has stimulated increasing research interest in learning algorithms for solving MILP problems. These works can be broadly divided into two categories: \emph{learning decision strategies within MILP solvers}, and \emph{learning primal heuristics}.

The first approach investigates the use of ML to learn to make decisions inside a MILP solver, which is typically built upon a general branch-and-bound framework. The learned policies can be either cheap approximations of existing expensive methods, or more sophisticated strategies that are to be discovered. Related works include: learning to select branching variables \cite{Khalil_LeBodic_Song_Nemhauser_Dilkina_2016, balcan2018learning, gasse2019exact}, learning to select branching nodes \cite{he2014learning}, learning to select cutting planes \cite{tang2020reinforcement}, and learning to optimize the usage of primal heuristics \cite{khalil2017learning, chmiela2021learning, hendel2022adaptive}.

The \emph{learning primal heuristics} approach is to learn algorithms to produce primal solutions for MILPs. Previous works in this area typically use ML methods to develop \emph{large neighborhood search} (LNS) heuristics. Within an LNS scheme, ML models are trained to predict ``promising" solution neighborhoods that are expected to contain high-quality solutions. In \cite{Ding_Zhang_Shen_Li_Wang_Xu_Song_2020}, the authors trained neural networks to directly predict solution values of binary variables, and then applied the LB heuristic to explore the solution neighborhoods around the predictions. The work \cite{nair2020solving} also uses neural networks to predict partial solutions. The subproblems defined by fixing the predicted partial solutions are solved by a MILP solver. The work \cite{sonnerat2021learning} proposes a LNS heuristic based on a ``learn to destroy" strategy, which frees part of the current solution. The variables to be freed are selected by trained neural networks using imitation learning. Note that their methods rely on parallel computation, which makes the outcome of the framework within a non-parallel environment less clear. In \cite{song2020general}, the authors proposed a decomposition-based LNS heuristic. They used imitation learning and reinforcement learning to decompose the set of integer variables into subsets of fixed size. Each subset defines a subproblem. The number of subsets is fixed as a hyperparameter.

Note that the learning-based LNS methods listed above directly operate on the integer variables, i.e., the predictions of ML models are at a variable-wise level, which still encounters the intrinsic combinatorial difficulty of the problem and limits their generalization performances on generic MILPs. Moreover, the learning of these heuristics is mostly based on the extraction of static features of the problem, the dynamic statistics of the heuristic behavior of the solver being barely explored. In our work, we aim at avoiding directly making predictions on variables. Instead, we propose to guide the (local) search by learning how to control the neighborhood size at an instance-wise level. To identify promising solution neighborhoods, our method exploits not only the static features of the problem, but also the dynamic features collected during the solution process as a sequential approach.

In the literature, there has also been an effort to learn algorithms for solving specific combinatorial optimization problems \cite{hottung2019neural, nazari2018reinforcement, kool2018attention, dai2017learning, bello2016neural, liu2021learning}. For a detailed overview of ``learn to optimize", see \cite{bengio2021machine}.


\section{Preliminaries}
\label{sec:background}
\subsection{Local branching}
\label{sed:lb}

We consider a MILP problem with 0–1 variables of the form

\begin{align}
 (P)~~~\min~~~ &\boldsymbol{c}^T\boldsymbol{x}\\
\text{s.t.}~~~  &\boldsymbol{Ax} \le \boldsymbol{b},\\
& x_j \in \{0,1\}, ~\forall j\in \mathcal{B},\\
& x_j \in \mathbb{Z}^+, ~\forall j\in \mathcal{G},~~~x_j \ge 0, ~\forall j\in \mathcal{C},
\end{align}
where the index set of decision variables $\mathcal{N}:=\{1,\ldots,n\}$ is partitioned into $\mathcal{B}, \mathcal{G}, \mathcal{C}$, which are the index sets of binary, general integer and continuous variables, respectively.

Note that we assume the existence of binary variables, as one of the basic building blocks of our method---namely, the local branching heuristic---is based on this assumption. However, this limitation can be relaxed and the local branching heuristic can be extended to deal with general integer variables, as proposed in \cite{bertacco2007general}.

Let $\bar\xx$ be a feasible \emph{incumbent} solution for $(P)$, and let $\mathcal{\overline{S}} = \{ j\in \mathcal{B} : \bar x_j=1 \}$ denote the binary support of $\bar\xx$. For a given positive integer parameter $k$, we define the neighborhood $N(\bar\xx, k)$ as the set of the feasible solutions of $(P)$ 
satisfying the \emph{local branching constraint}

\begin{equation}
\label{eq_phi}
\Delta(\boldsymbol{x}, \bar{\boldsymbol{x}}) = \sum_{j\in \mathcal{B}\setminus \mathcal{\overline S} }x_j + \sum_{j\in \mathcal{\overline S}}(1-x_j) \le k.
\end{equation}

In the relevant case in which solutions with a small binary support are considered (for example, in the famous traveling salesman problem only $n$ or the $O(n^2)$ variables take value 1), the asymmetric form of local branching constraint is suited, namely
\begin{equation}
\label{eq_phi2}
\Delta(\boldsymbol{x}, \bar\xx) = \sum_{j\in \mathcal{\overline S}}(1-x_j) \le k.
\end{equation}

The local branching constraint can be used in an exact branching scheme for $(P)$. Given the incumbent solution $\bar\xx$, the solution space with the current branching node can be partitioned by creating two child nodes as follows:
$$\text{Left: }~\Delta(\xx, \bar\xx) \le k,  ~~~~~\text{Right: } ~ \Delta(\xx, \bar\xx) \ge k + 1.$$
\subsection{The neighborhood size optimization problem}

For a \emph{neighborhood size} parameter $k \in \mathbb{Z}^{+}$, the LB algorithm obtained from choosing $k$ can be denoted as $\mathcal{A}_{k}$.  Given a MILP instance $\pp$, with its incumbent solution $\bar\xx$, the neighborhood size optimization problem over $k$ for one iteration of $\mathcal{A}_{k}$ is defined as
\begin{align}
     \min~~~ &\mathcal{C} ( \pp, \bar\xx ; \mathcal{A}_k )\\
    \text{s.t.}~~~  & k\in \mathbb{Z}^{+}
    ,
\end{align}
where $\mathcal{C} ( \pp, \bar\xx ; \mathcal{A}_k )$ measures the ``cost" of $\mathcal{A}_k $ on instance $(\pp, \bar\xx )$ as a trade-off between execution speed and solution quality.

In practice, a run of the LB algorithm consists of a sequence of LB iterations. To maximize the performance of the LB algorithm, a series of the above optimization problems need to be solved. Since the cost function $\mathcal{C}$ is unknown, those problems cannot be solved analytically. In general, the common strategy is to evaluate some trials of $k$ and select the most performing one. This is often done by using black-box optimization methods \cite{audet2017derivative}. As the evaluation of each setting involves a run of $\mathcal{A}_{k}$ and the best $k$ is instance-specific, those methods are not computationally efficient enough for online use. That is why the original LB algorithm initializes $k$ with a fixed small value and adapts it conservatively by a deterministic strategy.

Currently, learning from experiments and transferring the learned knowledge from solved instances to new instances is of increasing interest and somehow accessible. In the next section, we will introduce a new strategy for selecting the neighborhood size $k$ by using data-driven methods.
\section{Learning methods}
\label{sec:method}
Next, we present our framework for learning the neighborhood size in the LB scheme. The original LB algorithm chooses a conservative value for $k$ as default, with the aim of generating a easy-to-solve subproblem for general MILPs. However, as discussed in Section \ref{sec:intro}, our observation shows that the ``best" $k$ is dependent on the particular MILP instance. Hence, in order to optimize the performance of the LB heuristic, we aim at devising new strategies to learn how to tailor the neighborhood size for a specific instance. In particular, we investigate the dependencies between the state of the problem (defined by a set of both static and dynamic features collected during the LB procedure, e.g., context of the problem, incumbent solution, solving status, computation cost, etc.) and the size of the LB neighborhood. 

Our framework consists of a two-phase strategy. In the first phase, we define a regression task to learn the neighborhood size for the first LB iteration.  Within our method, this size is predicted by a pretrained regression model. 
For the second phase, we leverage reinforcement learning (RL) \cite{sutton1998introduction} and train a policy to dynamically adapt the neighborhood size at the subsequent LB iterations. The exploration of each LB neighborhood is the same as in the original LB framework, and a generic black-box MILP solver is used to update the incumbent solution. 
The overall scheme is exact in nature although turning it into a specialized heuristic framework is trivial (and generally preferrable).

\subsection{Scaled regression for local branching}
\label{sec:learning_k0}
For intermediate LB iterations, the statistics of previous iterations (e.g., value of the previous $k$, solving statistics, etc.) are available. One can then take advantage of this information and exploit the learning methods based on dynamic programming (e.g., reinforcement learning). Section \ref{sec:adapt_k} will address this case. However, for the first LB iteration, there is no historical information available as input. In this section, we will show how to define a regression task to learn the first $k$ from the context of the problem and the incumbent solution.

Let $\mathcal{S}$ denote the set of available features of the MILPs before the first LB iteration. We aim to train a regression model $f : \mathcal{S} \to \mathbb{R}$ that maps the features of a MILP instance $\ss$ to $k_0^*$, the label of best $k_0$. However, the label $k_0^*$ is unknown, and we do not have any existing method to compute the exact $k_0^*$. To generate labels, we first define a metric for assessing the performance of a LB algorithm $\mathcal{A}_k$, and then use black-box optimization methods to produce approximations of $k_0^*$ as labels.


\subsubsection{Approximation of the best $k_0$ }
To define a cost metric for $\mathcal{A}_k$, we consider two factors. The first factor is the computational effort (e.g., CPU time) to solve the sub-MILP defined by the LB neighborhood, while the second factor is the solution quality (e.g., the objective value of the best solution). To quantify the trade-off of speed and quality, the cost metric can be defined as
\begin{align}
c^{k_0} = \alpha * t^{k_0}_{scaled} + (1-\alpha) * o^{k_0}_{scaled}, \label{eq-9}
\end{align}
where $t^{k_0}_{scaled} \in [0, 1]$ is the scaled computing time for solving the sub-MILP, $o^{k_0}_{scaled} \in [0, 1]$ is the scaled objective of sub-MILP, and $\alpha$ is a constant.

\def\myhat{}  

Given $c^{k_0}$, the label $k_0^*$ can be defined as
\begin{align}
   k_{0}^{*} = \argmin_{ k_{0} \in \mathbb{Z}^{+}} c^{k_0},
\end{align}
and is usually evaluated through black-box optimization methods: Given the time limit 
and a collection of training instances of interest, one evaluates the LB algorithm introduced in Section \ref{sec:background} with different values of $k_0$,
and $k^*_0$ is estimated by choosing the value with the best performance assessed by (\ref{eq-9}), which is typically the largest $\myhat{k}_0$ such that the resulting sub-MILP can still be solved to optimality within the time limit. Since the evaluation process for each instance is quite expensive, we propose to approximate it through regression. 


\subsubsection{Regression for learning $k_0$}
With a collected dataset $\mathcal{D}=({\ss}_i, {\myhat{k}}^*_{0i})_{i=1}^N$ with $N$ instances, a regression task can be analyzed to learn a mapping from the state of the problem to the estimated $\myhat{k}_0^{*}$. The regression model $f_\theta (\ss)$ can be obtained by solving

\begin{align}
    \myhat{\theta^*} = \argmin_{\theta \in \Theta} \frac{1}{N} \sum_{i=1}^{N}{ \mathcal{L}(f_\theta (\ss_i), \myhat{k}^*_{0i}) },
\end{align}
where $\mathcal{L}(f_\theta (\ss_i), \myhat{k}^*_{0i})$ defines the loss function, a typical choice for regression task being the mean squared error.

\subsubsection{The scaled regression task}
Let $\xx'$ be the optimal linear programming (LP) fractional solution without local branching constraint, and let $k'$ be the value of the left-hand side of the local branching constraint evaluated using $\xx'$. Specifically, $k'$ is computed by
\begin{align}
    k' = \Delta(\xx', \bar{\xx}).
\end{align}

As discussed in \cite{fischetti2014proximity}, any $k \ge k'$ is likely to be useless as  the LP solution  after adding the LB constraint would be unchanged.  
Hence, $k'$ provides an upper bound for $k$. 

We can therefore parametrize $k$ as 
\begin{align}
    k = \phi \: k',
\end{align}
where $\phi \in (0,1)$. Now, we define the regression task over a scaled space $\phi \in (0,1)$ instead of directly over $k$.

Given $\myhat{k}^*_{0}$ and $k'_{0}$, the label is easily computed by
\begin{align}
    \myhat{\phi}^*_{0} = \frac{\myhat{k}^*_{0}}{k'_{0}}.
\end{align}

The regression problem reduces to
\begin{align}
    \myhat{\theta^*} = \argmin_{\theta \in \Theta} \frac{1}{N} \sum_{i=1}^{N}{ \mathcal{L}(f_\theta (\ss_i), \myhat{\phi}^*_{0i}) },
\end{align}
where $\mathcal{L}(f_\theta (\ss_i), \myhat{\phi}^*_{0i})$ defines the loss function.

\paragraph{MILP representation}
We represent the state $\ss$ as a bipartite graph $(\mathbf{C}, \mathbf{E}, \mathbf{V})$ \cite{gasse2019exact}. Given a MILP instance, let $n$ be the number of variables with $d$ features for each variable, $m$ be the number of constraints with $q$ features for each constraint. The variables of the MILP, with $\mathbf{V} \in \mathbb{R}^{n\times d}$ being their feature matrix, are represented on one side of the graph. On the other side are nodes corresponding to the constraints with $\mathbf{C} \in \mathbb{R}^{m\times q}$ being their feature matrix. A constraint node $i$ and a variable node $j$ are connected by an edge $(i,j)$ 
if variable $i$ appears in constraint $j$ in the MILP model.
Finally, $\mathbf{E} \in \mathbb{R}^{m\times n\times e}$ denotes the tensor of edge features, with $e$ being the number of features for each edge.

\begin{algorithm}[t]
\caption{LB with scaled regression}
\label{alg:the_alg_1}
\SetAlgoLined
\KwIn{instance dataset $\mathcal{P} = \{ {\pp}_i \}_{i=1}^M $}

 \For{instance ${\pp}_i \in \mathcal{P}$}{ 
    ~~0. initialize the state $\ss$ with an initial solution $\bar{\xx}$ \; \\
    1. solve the LP relaxation and get solution $\boldsymbol{x'}$ \; \\
    2. compute $k' = \Delta(\boldsymbol{x'}, \bar{\boldsymbol{x}})$ \; \\
    3. predict $\phi_0 = f_{\myhat{\theta}} (\ss) $ by the regression model \; \\
    4. compute $k_0 = \phi_0 \: k'$ \; \\
    5. apply $k_0$ to execute the first LB iteration\; \\
    6. update the incumbent $\bar{\xx}$ and continue LB algorithm with its default setting\;  \\
    \Repeat {termination condition is reached} {
        execute the next LB iteration
    }
    }
\end{algorithm}

\paragraph{Regression model}
Given that states are represented as graphs, with arbitrary size and topology, we propose to use \emph{graph neural networks} (GNNs) \cite{gori2005new, hamilton2017representation} to parameterize the regression model. Indeed, GNNs are size-and-order invariant to input data, i.e., they can process graphs of arbitrary size, and the ordering of the input elements is irrelevant. Another appealing property of GNNs is that they exploit the sparsity of the graph, which makes GNNs an efficient model for embedding MILP problems that are typically very sparse \cite{gasse2019exact}.

Our GNN architecture consists of three modules: the input module, the convolution module, and the output module. In the input layer, the state $\ss$ is fed into the GNN model. The input module embeds the features of the state $\ss$. The convolution module propagates the embedded features with graph convolution layers. In particular, our graph convolution layer applies the message passing operator, defined as

\begin{align}
    \mathbf{v}_i^{(h)} = f_{\theta}^{(h)} \left( \mathbf{v}_i^{(h-1)}, \sum_{j \in \mathcal{N}(i)} g_{\phi}^{(h)}\left(\mathbf{v}_i^{(h-1)}, \mathbf{v}_j^{(h-1)},\mathbf{e}_{j,i}\right) \right),
\end{align}
where $\mathbf{v}^{(h-1)}_i \in \mathbb{R}^d$ denotes the feature vector of node $i$ from layer $(h-1)$, $\mathbf{e}_{j,i} \in \mathbb{R}^m$ denotes the feature vector of edge $(j, i)$ from node $j$ to node $i$ of layer $(h-1)$, and $\mathbf{f}_{\theta}^{(h)}$ and $ \mathbf{g}_{\phi}^{(h)}$ denote the embedding functions in layer $h$.

For a bipartite graph, a convolution layer is decomposed into two half-layers: one half-layer propagates messages from variable nodes to constraint nodes through edges, and the other one propagates messages from constraint nodes to variable nodes.
The output module embeds the features extracted from the convolution module and then applies a pooling layer, which maps the graph representation into a single neuron. The output of this neuron is the prediction of $\phi_0$. 

\subsubsection{LB with scaled regression}
Our refined LB heuristic, \emph{LB with scaled regression}, is obtained when $k_0$ is predicted by the regression model. The pseudocode of the algorithm is outlined in Algorithm \ref{alg:the_alg_1}. 

\subsection{Reinforced neighborhood search}
\label{sec:adapt_k}
In this section, we leverage reinforcement learning (RL) to adapt the neighborhood size iteratively. We first formulate the problem as a Markov Decision Process (MDP) \cite{howard1960dynamic}.
Then, we propose to use policy gradient methods to train a policy model. 
\subsubsection{Markov Decision Process}

	Given a MILP instance with an initial feasible solution, the procedure can be formulated as a MDP, wherein at each step, a neighborhood size is selected by a policy model and applied to run a LB iteration. In principle, the state space $\mathcal{S}$ is the set of all the features of the MILP model and its solving statistics, which is combinatorial and arbitrarily large. To design an efficient RL framework for this problem, we choose a compact set of features from the solving statistics to construct the state. These features characterize the progress of the optimization process and are instance-independent, allowing broader generalization across instances.
	
	For the action space State($\mathcal{A}$), instead of directly selecting a new $k$, we choose to adapt the value of $k$ of the last LB iteration. The set of possible actions consists of four options
	\begin{align}
	    \{+k_{step},0, -k_{step}; reset \},
	\end{align}
	where ``$+k_{step}$" means increasing $k$ by $k = k + k_{step} * k$, ``$-k_{step}$" means decreasing $k$ by $k = k - k_{step}*k$, ``$0$" denotes keeping $k$ without any change, and ``$reset$" means resetting $k$ to a default value. The policy $\pi_k$ maps a state to one of the four actions. The step size $k_{step}$ is a hyperparameter of the algorithm.
	
	\begin{figure}[H]
          \centering
          \includegraphics[width=0.8\linewidth]{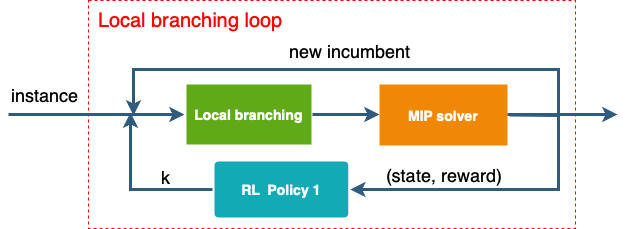}
          \caption{RL framework for adapting $k$ }
          \label{fig:plot_klb}
    \end{figure}
	
	The compact description of states and actions offers several advantages. First of all, it simplifies the MDP formulation and makes the learning task easier. In addition, it allows for the use of simpler function approximators, which is critical for speeding up the learning process. In Section \ref{sec:result}, we will show---by training a simple linear policy model using the off-the-shelf policy gradient method---that the resulting policy can significantly improve the performance of the LB algorithm.
	
	By applying the updated $k$, the next LB iteration is executed with time limit $t_{limit}$. Then, the solving sub-MILP statistics are collected to create the next state. In principle, the reward $r_k$ is formulated according to the outcome of the last LB iteration, e.g., the computing time and the quality of the incumbent solution. 
	
	To maximize the objective improvement and minimize the solution time of the LB algorithm, we define the combinatorial reward as
    \begin{equation}
        r_k = {o}_{imp} * (t_{max} - t_{elaps}),
    \end{equation}
    where ${o}_{imp}$ denotes the objective improvement obtained from the last LB iteration, $t_{max}$ is the global time limit of the LB algorithm, and $t_{elaps}$ is the cumulated running time.
    
    The definition above is just one possibility to build a MDP for the LB heuristic. Actually, defining a compact MDP formulation is critical for constructing efficient RL algorithms for this problem. 

\subsubsection{Learning strategy}
For training the policy model, we use the \emph{reinforce} policy gradient method \cite{sutton2000policy}, which allows a policy to be learned without any estimate of the value functions.

The refined LB heuristic, \emph{reinforced neighborhood search for adapting $k$} is obtained when the neigborhood size $k$, is dynamically adapted by the RL policy. The pseudocode of the algorithm is outlined in Algorithm \ref{alg:the_alg_2}.

\begin{algorithm}[t]
\caption{Reinforced neighborhood search for adapting $k$}
\label{alg:the_alg_2}
\SetAlgoLined
\KwIn{instance dataset $\mathcal{P} = \{ {\pp}_i \}_{i=1}^M $}

 \For{instance ${\pp}_i \in \mathcal{P}$}{ 
    ~~0. initialize the state $\ss$ with an initial solution $\bar{\xx}$ \; \\
    1. compute $k_0$ by the procedure in Algorithm 1 or set $k_0$ by a default value \; \\
    2. apply $k_0$ to execute the first LB iteration\; \\
    3. collect the new state $\ss$ and the incumbent  $\bar{\xx}$\;  \\
    \Repeat {termination condition is reached} {
        ~~update $k$ by policy $\pi_k (\ss)$ \; \\
        apply $k$ and execute the next LB iteration \; \\
        collect the new state $\ss$ with the incumbent $\bar{\xx}$\;
    }
    }
\end{algorithm}

\subsection{Further improvement by adapting LB node time limit}
\label{sec:adapt_t}
\begin{figure}[H]
          \centering
          \includegraphics[width=0.8\linewidth]{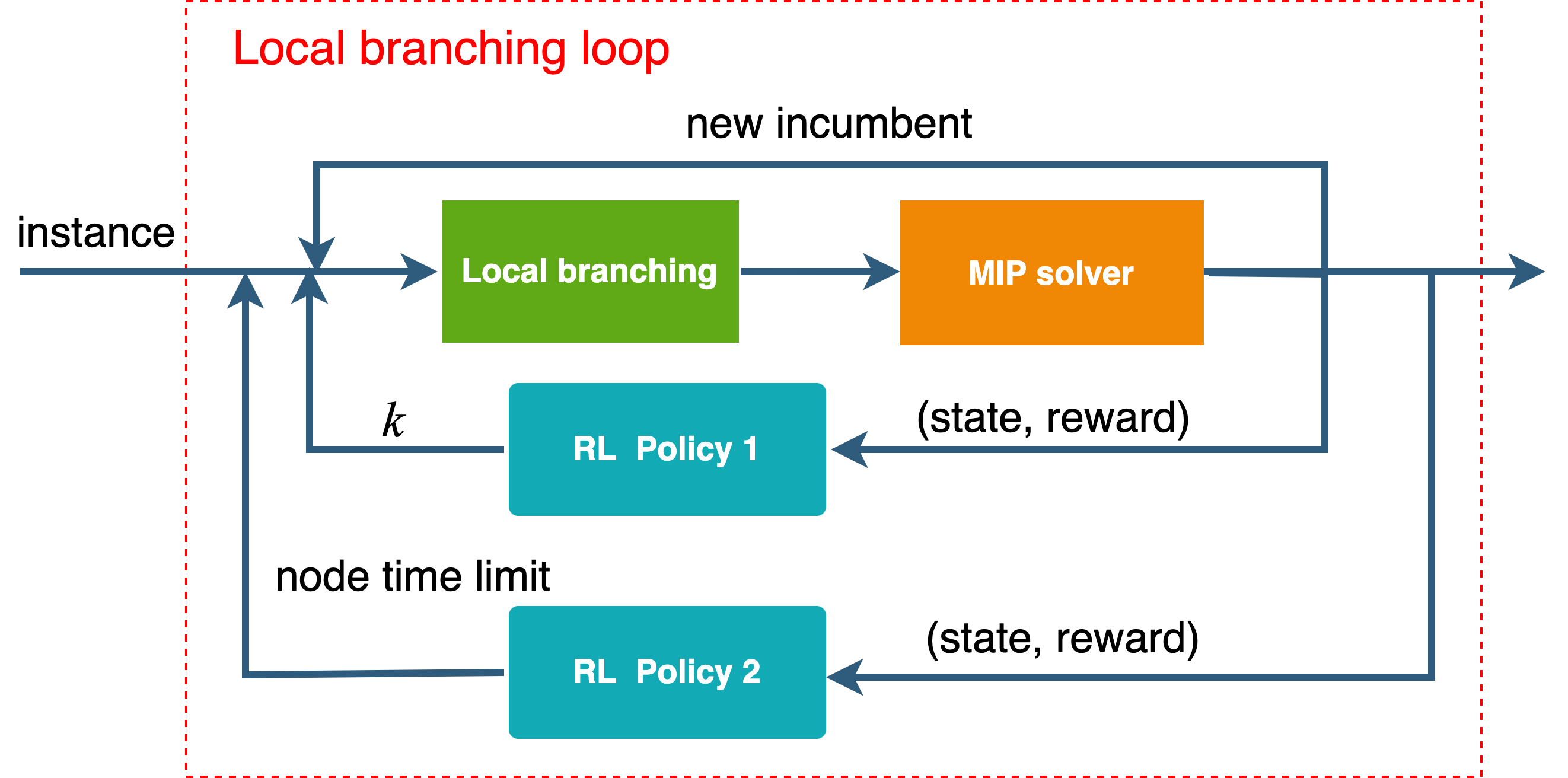}
          \caption{RL framework for adapting the time limit for solving the LB subproblem.}
          \label{fig:plot_tlb}
    \end{figure}

In the previous section, we suppose that the time limit for solving each local branching neighborhood is given. Indeed, the setting of time limit for each ``LB node" decides how much computational effort is spent on each LB node, and therefore customizes how the overall time limit is split. E.g., a larger ``LB node" time limit allows more computation for solving that LB node, which potentially leads to a larger $k$. However, it should be noted that the relation between the value of ``LB node" time limit and the actual performance of LB scheme is not known. Therefore, we seek to address this issue by leveraging reinforcement learning again to learn policies for tailoring the time limit for each LB node. This new ``RL-t" loop is built on top of the inner ``RL-k" loop for adapting the LB neighborhood size. The scheme is shown in Figure \ref{fig:plot_tlb}.

To design the action space for adapting the time limit of each LB node, we apply a similar setting as used for adapting $k$. The set of possible actions consists of four options
	\begin{align}
	    \{*t_{step},0, /t_{step}; reset \},
	\end{align}
	where ``$*t_{step}$" means increasing $t$ by $t = t * t_{step}$, ``$/t_{step}$" means decreasing $t$ by $t = t / t_{step}$, ``$0$" denotes keeping $t$ without any change, and ``$reset$" means resetting $t$ to a default value. The policy $\pi_t$ maps a state to one of the four actions.
	
For the reward design of the ``RL-t" loop, the reward also takes into account both the cost of computation spent on solving the LB subproblem and the solution quality. Specifically, the reward $r_t$ for the outer RL-t loop consists of two components. The first component $r_{1}$ maintains the same reward $r_k$ as used in the ``RL-k" loop, which attempts to maximize the improvement of the objective by the last LB iteration while minimizing the computing time. In addition, the LB environment explicitly returns a reward signal as a second component to penalize those LB iterations where the subproblem is too large to solve and there is no objective improvement within the time limit.
The penalty reward $r_p$ is a binary signal defined as
\begin{align*}
        r_p= 
        \left\{
        \begin{array}{rl}
            1 & \ \ \ \text{if the LB subproblem is not solved and no improving} 
            \\
              &  \ \ \ \text{solution is returned,}
            \\
            0 & \ \ \ \text{otherwise.}
        \end{array}
        \right.
\end{align*}

Formally, the reward $r_t$ is a combination of two components, defined as
\begin{align*}
    r_t = \beta_1 r_1 + \beta_2 r_2 \;
\end{align*}
where $r_1 = r_k$, $r_2=r_p$, $\beta_1 > 0$, $\beta_2 > 0$. 
    
For training the policy $\pi_t (\ss)$, we use the same RL method as used for training the policy $\pi_k (\ss)$. Note that we fix the pretrained $\pi_k (\ss)$ policy while tuning $\pi_t (\ss)$ to make the learning process more stable.

A new LB alogorithm, \emph{hybrid reinforced neighborhood search}, is obtained when the neigborhood size $k$ and the time limit $t$ for each LB node are dynamically adapted by the RL policies. The pseudocode of the algorithm is outlined in Algorithm \ref{alg:the_alg_3}.

\begin{algorithm}[t]
\caption{Hybrid reinforced neighborhood search for adapting $k$ and $t$}
\label{alg:the_alg_3}
\SetAlgoLined
\KwIn{instance dataset $\mathcal{P} = \{ {\pp}_i \}_{i=1}^M $}

 \For{instance ${\pp}_i \in \mathcal{P}$}{ 
    ~~0. initialize the state $\ss$ with an initial solution $\bar{\xx}$ \; \\
    1. compute $k_0$ by the procedure in Algorithm 1 or set $k_0$ by a default value \; \\
    2. apply $k_0$ to execute the first LB iteration\; \\
    3. collect the new state $\ss$ and the incumbent  $\bar{\xx}$\;  \\
    \Repeat {termination condition is reached} {
        ~~update $k$ by policy $\pi_k (\ss)$ \; \\
        update $t$ by policy $\pi_t (\ss)$ \; \\
        apply $k$, $t$ and execute the next LB iteration \; \\
        collect the new state $\ss$ with the incumbent $\bar{\xx}$\;
    }
    }
\end{algorithm}
\section{Experiments}
\label{sec:result}
In this section, we present the details of our experimental results over five MILP benchmarks. 
We compare different settings of our approach against the original LB algorithm, using SCIP \cite{GamrathEtal2020ZR} as the underlying MILP solver.

\subsection{Data collection}

\subsubsection{MILP instances}
We evaluate on five MILP benchmarks: set covering (SC) \cite{balas1980set}, maximum independent set (MIS) \cite{bergman2016decision}, combinatorial auction (CA) \cite{leyton2000towards}, generalized independent set problem (GISP) \cite{hochbaum1997forest, colombi2017generalized}, and MIPLIB 2017 \cite{gleixner2021miplib}. The first three benchmarks are used for both training and evaluation. For SC, we use instances with $5000$ rows and $2000$ columns. For MIS, we use instances on Barabási–Albert random graphs with $1000$ nodes. For CA, we use instances with 4000 items and 2000 bids. In addition, to evaluate the generalization performance on larger instances, we also use a larger dataset of instances with doubled size for each benchmark, denoted by LCA, LMIS, LCA. The larger datasets are only used for evaluation. 

For GISP, we use the public dataset from \cite{chmiela2021learning}. For MIPLIB, we select binary integer linear programming problems from MIPLIB 2017. Instances from GISP and MIPLIB are only used for evaluation.

For each instance, an initial feasible solution is required to run the LB heuristic. We use two initial incumbent solutions: (1) the first solution found by SCIP; (2) an intermediate solution found by SCIP, typically the best solution obtained by SCIP at the end of the root node computation, i.e., before branching.


\subsubsection{Data Collection for regression}
To collect data for the scaled regression task, one can use black-box optimization methods to produce the label $\phi_0^{*}$. As the search space has only one dimension, we choose to use the grid search method. In particular, given a MILP instance, an initial incumbent $\bar\xx$, the LP solution $\xx'$,  and a time limit for a LB iteration, we evaluate $\phi_{0}$ from $(0, 1)$ with a resolution limit $0.01$. For each $\phi_{0}$, we compute the actual neighborhood size by $k_0 = k'\: \phi_0$, where $k' = \Delta(\xx', \bar\xx)$. Then, $k_0$ is applied to execute an iteration of LB. From all the evaluated $\phi_{0}$, the one with best performance is chosen as a label $\phi_{0}^{*}$.




The state $\mathbf{\ss}$ consists of context features of the MILP model and the incumbent solution. The state $\ss$ together with the label $k^*$ construct a valid data point ($\ss$, $k^*$).

\subsection{Experimental setup}

\subsubsection{Datasets}
For each reference set of SC, MIS and CA problems, we generate a dataset of $200$ instances, and split each dataset into training ($70\%$), validation ($10\%$), and test ($20\%$) sets. For larger instances, we generate $40$ instances of LSC, LMIS and LCA problems, separately. The GISP dataset contains $35$ instances. For MIPLIB, we select $29$ binary MILPs that are also evaluated by the original LB heuristic \cite{fischetti2003local}.

\subsubsection{Model architecture and feature design}

For the regression task, we apply the GNNs described in the paper with three modules. For the input module, we apply 2-layer perceptron with the \emph{rectified linear unit} (ReLU) activation function to embed the features of nodes. For the convolution module, we use two half-layers, one from nodes of variables to nodes of constraints, and the other one from nodes of constraints to nodes of variables. For the output module, we also apply 2-layer perceptron with the ReLU activation function. The pooling layer uses the \emph{sigmoid} activation function. All the hidden layers have $64$ neurons. 

Given the input $x\; {\in}\; \mathbb{R}$, the $\Sigmoid$ and $\ReLU$ functions are defined as
\begin{align}
       \Sigmoid (x) &= \frac{\exp(x)}{\exp(x) + 1}\;, \\
       \ReLU (x) &= \max(0,x)\;.
   \end{align}

For the bipartite graph representation, we reference the model used in \cite{gasse2019exact}. The features in the bipartite graph are listed in Table \ref{table:bipartite-features}. In practice, if the instance is a pure binary MILP, one can choose a more compact set of features to accelerate the training process--for example, the features describing the type and bound of the variables can be removed.

\begin{table}[htbp!]
    \vspace{10pt}
    \centering
    \begin{tabular}{c l l}
    \multicolumn{1}{c}{Tensor} & \multicolumn{1}{l}{Feature} & \multicolumn{1}{l}{Description} \\
    \toprule
    \multirow{1}{*}{$\mathbf{C}$} & bias & Bias value, normalized with constraint coefficients. \\
    \midrule
    \multirow{1}{*}{$\mathbf{E}$}
    & coef & Constraint coefficient, normalized per constraint. \\
    \midrule
    \multirow{10}{*}{$\mathbf{V}$}
    & coef & Objective coefficient, normalized. \\
    \cmidrule{2-3}
    & binary & Binary type binary indicator. \\
    \cmidrule{2-3}
    & integer & Integer type indicator. \\
    \cmidrule{2-3}
    & imp\_integer & Implicit integer type indicator \\
    \cmidrule{2-3}
    & continuous & Continuous type indicator.  \\
    \cmidrule{2-3}
    & has\_lb & Lower bound indicator. \\
    \cmidrule{2-3}
    & has\_ub & Upper bound indicator. \\
    \cmidrule{2-3}
    & lb & Lower bound. \\
    \cmidrule{2-3}
    & ub & Upper bound. \\
    \cmidrule{2-3}
    & sol\_val & Solution value.  \\
    \bottomrule
    \end{tabular}
    \caption{Description of the features in the bipartite graph $\mathbf{\ss} = (\mathbf{C}, \mathbf{E}, \mathbf{V})$.}
    \label{table:bipartite-features}
\end{table}

For the two RL policies, we apply a linear model with seven inputs and four outputs. The set of features used by the RL policies is listed in Table \ref{table:rl-features}.

\begin{table}[htbp!]
    \vspace{10pt}
    \centering
    \begin{tabular}{l l}
    \multicolumn{1}{l}{Feature} & \multicolumn{1}{l}{Description} \\
    \toprule
    optimal & Sub-MILP is solved and the incumbent is updated, indicator. \\
    \midrule
    infeasible & Sub-MILP is proven infeasible, indicator. \\
    \midrule
    improved & Sub-MILP is not solved but the incumbent is updated, indicator.  \\
    \midrule
    not\_improved& Sub-MILP is not solved and no solution is found, indicator. \\
    \midrule
    diverse & No improved solution is found for two iterations, indicator. \\
    \midrule
    t\_available & time available before the time limit of the sub-MILP is reached \\
    \midrule
    obj\_improve & improvement of objective.  \\
    \midrule
    \end{tabular}
    \caption{Description of the input features of the RL policy.}
    \label{table:rl-features}
\end{table}

\subsubsection{Training and evaluation}
For the regression task, the model learns from the features of the MILP formulation and the incumbent solution. We use the mean squared error as the loss function.
We train the regression model with two scenarios: the first one trains the model on the training set of SC, MIS and CA separately, the other one trains a single model on a mixed dataset of the three training sets. The models trained from the two scenarios are compared on the three test sets. 

For the RL task, since we only use the instance-independent features selected from solving statistics, the RL policy is only trained on the training set of SC, and evaluated on all the test sets.

To further evaluate the generalization performance with respect to the instance size and the instance type, the RL policies (trained on the SC dataset) and the regression model (trained on the SC, MIS and LCS datasets) were evaluated on GISP and MIPLIB datasets.


\subsubsection{Evaluation metrics}
We use two measures to compare the performance of different heuristic algorithms. The first indicator is the \emph{primal gap}. Let $\tilde{\xx}$ be a feasible solution , and $\tilde{\xx}_{opt}$ be the optimal (or best known) solution. The \emph{primal gap} (in percentage) is defined as
\begin{align*}
        \gamma (\tilde{\xx}) = 
       \frac{|c^T \tilde{\xx}_{\text{opt}} - c^T \tilde{\xx} |}{|c^T \tilde{\xx}_{\text{opt}}|} \times 100 
       ,
\end{align*}
where we assume the denominator is nonzero.    

For the second measure, we use the \emph{primal integral} \cite{berthold2013measuring}, which takes into account both the quality of solutions and the solving time required to find them. To define the primal integral, we first consider a \emph{primal gap function} $p(t)$ as a function of time, defined as
\begin{align*}
        p(t) = 
        \left\{
        \begin{array}{ll}
            1, & \ \ \ \text{if no incumbent until time $t$},\\
            \bar{\gamma} (\tilde{\xx}(t)), &\ \ \ \text{otherwise},
        \end{array}
        \right.
\end{align*}
where $\tilde{\xx}(t)$ is the incumbent solution at time $t$, and $\bar{\gamma} (\cdot) \in [0,1]$ is the \emph{scaled primal gap} defined by
\begin{align*}
        \bar{\gamma} (\tilde{x}) = 
        \left\{
        \begin{array}{ll}
            0, & \ \ \ \text{if $c^T \tilde{\xx}_{\text{opt}}=c^T \tilde{\xx}=0$},\\
            1, & \ \ \ \text{if $c^T \tilde{\xx}_{\text{opt}} \cdot c^T \tilde{\xx}<0$},\\
            \frac{|c^T \tilde{\xx}_{\text{opt}} - c^T \tilde{\xx} |}{\text{max}\{|c^T \tilde{\xx}_{\text{opt}}|, \ |c^T \tilde{\xx}| \}}, & \ \ \ \text{otherwise}.
        \end{array}
        \right.
\end{align*}
Let $t_{\text{max}} > 0$ be the time limit. The primal integral of a run is then defined as
\begin{align*}
        P(t_{\text{max}}) = 
        \int_{0}^{t_{\text{max}}} p(t) \, dt.
\end{align*}

\subsubsection{Details of experimental settings}
To tune the learning rate for training the regression model, we have experimented different learning rates from $10^{-5}$ to $10^{-1}$ and have chosen a learning rate of $10^{-4}$. We trained the model with a limit of 300 epochs.

For training the RL policies, we used the same method to tune the learning rate and have chosen a learning rate of $10^{-2}$ for $\pi_k$, and $10^{-1}$ for $\pi_t$ . We trained the RL policies with a limit of $300$ epochs.

For the hyperparameters, we have chosen $\alpha = 0.5$ as a trade-off between speed and quality for the cost metric for $k_0$. We used $
\beta_1 = \beta_2 = 1$ for the two components of $r_t$. We used $k_{step}=0.5$ for the action design in Section \ref{sec:adapt_k} and $t_{step}=2$ for the action design in Section \ref{sec:adapt_t}. We set a time limit of $10$ seconds for each LB iteration for all the compared algorithms.

Our code is written in Python $3.7$ and we use Pytorch $1.60$ \cite{paszke2019pytorch}, Pytorch Geometric $1.7.0$ \cite{fey2019fast}, PySCIPOpt $3.1.1$ \cite{MaherMiltenbergerPedrosoRehfeldtSchwarzSerrano2016}, SCIP $7.01$ \cite{GamrathEtal2020ZR} for developing our models and sovling MILPs.

\subsection{Results}
In order to validate our approach, we first implement LB as a heuristic search strategy for improving a certain incumbent solution using the MILP solver as a black-box. 
\subsubsection{Local branching search with adapting $k$}

We 
perform the evaluations of our framework on the following four settings:
\begin{itemize}
	\item \emph{lb-sr}: Algorithm \ref{alg:the_alg_1} with regression model trained by a homogenous dataset of \emph{SC}, \emph{MIS}, \emph{CA}, separately;
	\item \emph{lb-srm}: Algorithm \ref{alg:the_alg_1} with regression model trained by a mixed dataset of \emph{SC}, \emph{MIS}, \emph{CA};
	\item \emph{lb-rl}: Algorithm \ref{alg:the_alg_2} with setting $k_0$ by a default value;
	\item \emph{lb-srmrl}: combined algorithm using regression from Algorithm \ref{alg:the_alg_1} (with regression model trained by mixed dataset of \emph{SC}, \emph{MIS}, \emph{CA}) and RL from Algorithm \ref{alg:the_alg_2}.
\end{itemize}
We use the original local branching algorithm as the baseline. All the algorithms use SCIP as the underlying MILP solver and try to improve the initial incumbent with a time limit of $60s$. Our code and more details of the experiment environment are publicly available\footnote{\url{https://github.com/pandat8/ML4LB}}. 

The evaluation results for the basic SC, MIS, CA datasets are shown in Table \ref{tab:pi_basic} and Table \ref{tab:pg_basic}. Our first observation is that the learning based algorithms of our framework significantly outperform the original LB algorithm. Both the primal integral and the final primal gap of the four LB variants are smaller than those of the baseline over most datasets, showing improved heuristic behavior.
Note that, although the regression model trained by using supervised learning and the policy model trained by RL can be used independently, they benefit from being used together. As a matter of fact, the hybrid algorithm \emph{lb-srmrl} combining both methods achieves a solid further improvement and outperforms the other algorithms for most cases.

We also evaluate the impact of the choice of training set for the regression model. By comparing \emph{lb-sr} and \emph{lb-srm}, we observe that the regression model trained on a mixed dataset of \emph{SC}, \emph{MIS}, \emph{CA} exhibit a performance very close to that of the model trained on a homogeneous dataset. Indeed, the GNN networks we used embed the features of the MILP problem and its incumbent solution. In particular, one significant difference between our method and those of previous works is that, instead of training a separate model for each class of instances, our method is able to train a single model yielding competitive generalization performances across instances. This is because our models predict the neighborhood size at a instance-wise level, rather than making predictions on variables.

\begin{table}[t]
 \vspace{-2mm}
 \centering
 \footnotesize
 \tabcolsep=0.11cm
 \begin{tabular}{@{}lrrrrrrrrrrrr@{}}
 \toprule
          & \multicolumn{2}{c}{SC} &  \multicolumn{2}{c}{MIS}  &    \multicolumn{2}{c}{CA}  \\ \cmidrule(lr){2-3} \cmidrule(lr){4-5} \cmidrule(lr){6-7}          
Algo.       &   first       &   root	&	first	&	root    &	first	&	root	\\	\midrule
lb-base	    &	57.013 	    &	5.697	&	5.713	&	4.893	&	9.296	&	3.308\\	
lb-sr       &	3.643	    &	1.751	&	1.537	&	3.499	&	6.358	&	2.025\\	
lb-srm      &	4.338	    &	1.818	&	1.419	&	3.654	&	6.556	&	2.120\\	
lb-rl       &   17.304      &	4.513   &	2.616   &	2.616   &	4.450   &	2.042\\	
lb-srmrl    &	\B{2.991}   &\B{1.567}  &\B{1.244}  &\B{2.452}  &\B{3.202}  &\B{1.454}\\
 \bottomrule   
 \end{tabular}
 \caption{Primal integral (geomeric means) for SC, MIS, CA problems.}
 \label{tab:pi_basic}
\end{table}

\begin{table}[t]
 \vspace{-2mm}
 \centering
 \footnotesize
 \tabcolsep=0.11cm
 \begin{tabular}{@{}lrrrrrrrrrrrr@{}}
 \toprule
          & \multicolumn{2}{c}{SC} &  \multicolumn{2}{c}{MIS}  &    \multicolumn{2}{c}{CA}  \\ \cmidrule(lr){2-3} \cmidrule(lr){4-5} \cmidrule(lr){6-7}          
Algo.       &   first       &   root	&	first	    &	root    &	first	    &	root	\\	\midrule
lb-base	    &	1411.624     &	1.390    &	0.326	    &\B{0.193}    &	3.307	    &	1.218\\	
lb-sr       &	1.425	    &	1.102    &	\B{0.220}	& 0.210    &	2.460	    &	0.877\\	
lb-srm      &	1.645	    &	1.118    &	0.247	    & 0.232    &	2.466	    &	0.971\\	
lb-rl       &	6.876   	 &	1.621    &	0.449       & 0.263    &	0.558  	    &	\B{0.229} \\	
lb-srmrl    &\B{1.335}       &\B{0.550}   &	0.275       & 0.242    &\B{0.465}       &	0.311 \\
 \bottomrule   
 \end{tabular}
 \caption{Final primal gap (geometric means in percentage) for SC, MIS, CA problems.}
 \label{tab:pg_basic}
\end{table}

\emph{Broader Generalization}
Next, we evaluate the generalization performance with respect to the size and the type of instances. Recall that our regression model is trained on a randomly mixed dataset of SC, MIS and CA problems, and the RL policy is only trained on the training set of SC. We evaluate the trained models on larger instances (LSC, LMIS, LCA) and new MILP problems (GISP, MIPLIB). The results of evaluation on larger instances are shown in Table \ref{tab:pi_large} and Table \ref{tab:pg_large}, whereas the results on GISP and MIPLIB datasets, are shown in Table \ref{tab:pi_new} and Table \ref{tab:pg_new}. 

\begin{table}[t]
 \vspace{-2mm}
 \centering
 \footnotesize
 \tabcolsep=0.11cm
 \begin{tabular}{@{}lrrrrrrrrrrrr@{}}
 \toprule
          & \multicolumn{2}{c}{LSC} &  \multicolumn{2}{c}{LMIS}  &    \multicolumn{2}{c}{LCA}  \\ \cmidrule(lr){2-3} \cmidrule(lr){4-5} \cmidrule(lr){6-7}          
Algo.       &   first   &   root	&	first	&	root    &	first	&	root	\\	\midrule
lb-base	    &	59.303 	&	19.147	&	26.230	&	27.228	&	30.610	&	13.187\\	
lb-srm      &	3.642	&	2.350	&	4.746	&	9.175	&	19.639	&	6.869\\	
lb-rl       &	45.658  &	7.889   &	5.476   &	7.929   &	17.338  &	7.081\\	
lb-srmrl    &\B{2.790}  &\B{1.876}  &\B{1.819}  &\B{5.807}  &\B{10.078}   &\B{4.311}\\
  \bottomrule   
  \end{tabular}
 \caption{Primal integral (geometric means) for LSC, LMIS, LCA problems.}
  \label{tab:pi_large}
\end{table}

\begin{table}[htbp!]
 \vspace{-2mm}
 \centering
 \footnotesize
 \tabcolsep=0.11cm
 \begin{tabular}{@{}lrrrrrrrrrrrr@{}}
 \toprule
          & \multicolumn{2}{c}{LSC} &  \multicolumn{2}{c}{LMIS}  &    \multicolumn{2}{c}{LCA}  \\ \cmidrule(lr){2-3} \cmidrule(lr){4-5} \cmidrule(lr){6-7}          
Algo.       &   first       &   root	&	first	    &	root    &	first	    &	root	\\	\midrule
lb-base	    &	8105.973    &	4.221   &	18.633	    &	12.342  &	26.429	    &	5.156\\	
lb-srm      &	2.379	    &	1.242   &	0.230	    &	0.659   &	15.882	    &	1.961\\	
lb-rl       &	136.401     &	4.216   &	0.362       &  0.206    &	3.341       &	0.987   \\	
lb-srmrl    &\B{1.326}      &\B{0.777}  &\B{0.195}      &\B{0.205}  &\B{2.007}      &\B{0.152}    \\
  \bottomrule   
  \end{tabular}
 \caption{Final primal gap (geometric means in percentage) for LSC, LMIS, LCA problems.}
 \label{tab:pg_large}
\end{table}

\begin{table}[htbp!]
 \vspace{-2mm}
 \centering
 \footnotesize
 \begin{tabular}{@{}lrrrrrrrr@{}}
 \toprule
          & \multicolumn{2}{c}{GISP} &  \multicolumn{2}{c}{MIPLIB}  \\ \cmidrule(lr){2-3} \cmidrule(lr){4-5}      
Algo.       &   first   &   root	&	first	&	root    	\\	\midrule
lb-base	    &	19.833 	&	16.692	&	12.318	&	9.003	\\	
lb-srm      &\B{13.359}	&	10.173	&   11.319	&\B{6.863}	\\	
lb-rl       &	18.949  &	15.429  &	12.676  &	8.318    \\	
lb-srmrl    &   13.739  & \B{9.865} &\B{10.172} &   7.151  \\
  \bottomrule   
  \end{tabular}
 \caption{Primal integral (geometric means) for GISP and MIPLIB problems.}
 \label{tab:pi_new}
\end{table}

\begin{table}[htbp!]
 \vspace{-2mm}
 \centering
 \footnotesize
 \begin{tabular}{@{}lrrrrrrrr@{}}
 \toprule
          & \multicolumn{2}{c}{GISP} &  \multicolumn{2}{c}{MIPLIB}  \\ \cmidrule(lr){2-3} \cmidrule(lr){4-5}      
Algo.       &   first       &   root	    &	first	    &	root    	\\	\midrule
lb-base	    &	22.244 	&	17.171	 &	44.794	&	11.893	\\	
lb-srm      &\B{10.307}	&	7.704	 &	32.061 	&\B{8.730}	\\	
lb-rl       &	20.030  &	14.402   &	42.676  &	11.869  \\	
lb-srmrl    &   11.641  & \B{5.929}  &\B{20.575}  & 9.696   \\
  \bottomrule   
  \end{tabular}
 \caption{Final primal gap (geometric means in percentage) for GISP and MIPLIB problems.}
 \label{tab:pg_new}
\end{table}



Overall, all of our learning-based LB algorithms outperform the baseline, and the hybrid algorithm \emph{lb-srmrl} achieves the best performance on most datasets. These results show that our models, trained on smaller instances, generalize well both with respect to the instance size and, remarkably, across instances.

\subsubsection{Local branching search with adapting both k and t}
In this section, we analyze the results of our approach on adapting both the neighborhood size $k$ and time limit $t$ for each LB node. Again, local branching is implemented as an independent heuristic search scheme for improving a certain incumbent solution using a black-box MILP solver. 

We compare the \emph{lb-base} and our best algorithm \emph{lb-srmrl} outlined in the previous section with a new algorithm:
\begin{itemize}
	\item \emph{lb-srmrl-adapt-t}: Combined RL algorithm using regression from Algorithm \ref{alg:the_alg_1} (with regression model trained by mixed dataset of \emph{SC}, \emph{MIS}, \emph{CA}) and Hybrid RL from Algorithm \ref{alg:the_alg_3}.
\end{itemize}

In order to validate the effectiveness of our strategy for adapting $t$, we conducted the experiments on the two difficult MIP datasets (GISP, MIPLIB) with a longer global time limit, $600s$.  

\begin{figure}[htbp!]
    \centering
    
    \includegraphics[width=0.48\textwidth]{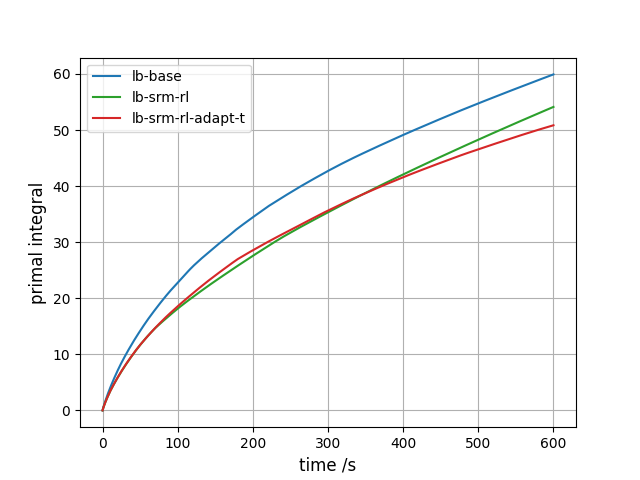}
    \hspace{0.02\textwidth}
    \includegraphics[width=0.48\textwidth]{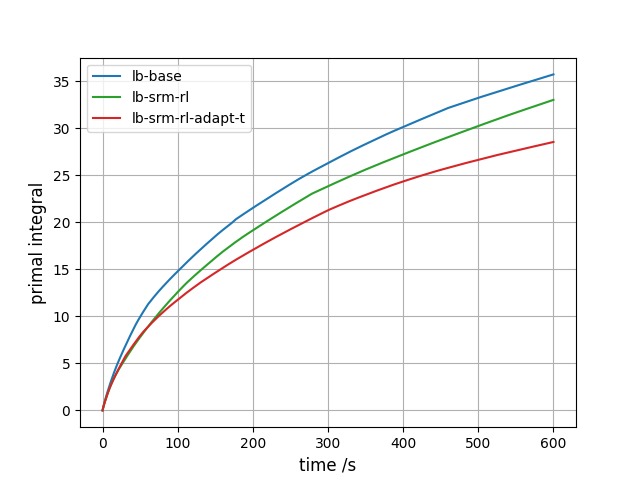}
    \caption{Evolution of the \emph{primal integral} (geometric means) over time on binary MIPLIB dataset. Left / right: using the first / root solution to start LB.}
    \label{fig:pi_new_long}
\end{figure}

The results are shown in Table \ref{tab:pi_new_long} and Table \ref{tab:pg_new_long}. In order to demonstrate the outcome of our LB algorithms on the solving progress of the instances over the running time, we plot the evolution of the average \emph{primal integral} on the MIPLIB dataset in Figure \ref{fig:pi_new_long}. The plot on the left reports the results of the instances initialized by the first solution found by SCIP, whereas the plot on the right reports the results of the same instances initialized by the best solution obtained at the end of the root node of the B\&B tree.

\begin{table}[htbp!]
 \vspace{-2mm}
 \centering
 \footnotesize
 \begin{tabular}{@{}lrrrrrrrr@{}}
 \toprule
          & \multicolumn{2}{c}{GISP} &  \multicolumn{2}{c}{MIPLIB}  \\ \cmidrule(lr){2-3} \cmidrule(lr){4-5}      
Algo.       &   first   &   root	&	first	&	root    	\\	\midrule
lb-base	    &	116.882 &	101.013	&	59.916	&	35.747	\\	
lb-srm-rl    &   108.345 &   91.600  &	54.116  &    33.038  \\
lb-srm-rl-adapt-t&\B{81.888}  &\B{68.042} &\B{50.843}&\B{28.558}    \\	
 \bottomrule   
 \end{tabular}
 \caption{Primal integral (geometric means) for GISP and MIPLIB problems with a time limit of $600s$ for each instance.}
 \label{tab:pi_new_long}
  
\end{table}

\begin{table}[htbp!]
 \vspace{-2mm}
 \centering
 \footnotesize
 \begin{tabular}{@{}lrrrrrrrr@{}}
 \toprule
          & \multicolumn{2}{c}{GISP} &  \multicolumn{2}{c}{MIPLIB}  \\ \cmidrule(lr){2-3} \cmidrule(lr){4-5}      
Algo.       &   first       &   root	    &	first	    &	root    	\\	\midrule
lb-base	    &	12.440 	&	11.039	&	36.860	&	\B{3.389}\\	
lb-srm-rl    &  12.230  &  11.991   &   20.532  &	5.273   \\
lb-srm-rl-adapt-t&\B{5.880}& \B{5.328} & \B{16.723}  &	4.396    \\	
 \bottomrule   
 \end{tabular}
 \caption{Final primal gap (geometric means in percentage) for GISP and MIPLIB problems with a time limit of $600s$ for each instance.}
 \label{tab:pg_new_long}
\end{table}

\newpage

From these results, we observe that: 1) All our learning-based LB algorithms converge faster than the LB baseline, showing improved heuristic behavior; 2) With our new strategy for adapting the LB node time limit $t$, the hybrid RL algorithm further improves its baseline with only a single RL policy for adapting $k$. The improvement becomes more and more significant as the solving time increases.


\section{Local branching as a primal heuristic within a MILP solver}
\label{sec:implementations}

Local branching can also be implemented as a refinement heuristic within a generic MILP solver. In this section, we present two possible implementations of the LB algorithms outlined in the previous section, to be used as a primal heuristic. One major difference between these implementations and those of the previous section is how the global MILP search is structured. Instead of applying LB as a metaheuristic strategy, we use LB within the MILP solver to improve the incumbent at certain nodes of B\&B tree. In particular, we considered the following two possible implementations:

\begin{itemize}
    \item executing LB primal heuristic only at the root node: the MILP solver calls local branching only at the root node or at the node where the first incumbent solution is found;
    \item executing LB primal heuristic at multiple nodes: the MILP solver checks if there is a new incumbent or not for every $f$ nodes in the B\&B tree, and calls local branching if that is the case. 
\end{itemize}

To configure the frequency $f$ of executing the LB primal heuristic for the second implementation, we have conducted a simple hyperparameter search for $f$ from the set of $\{1, 10, 100, 1000\}$, and the results showed that $100$ performed the best. Thus, we set $f=100$.

\begin{figure}[htbp!]
    \centering
    \hspace{0.02\textwidth}
    \includegraphics[width=0.60\textwidth]{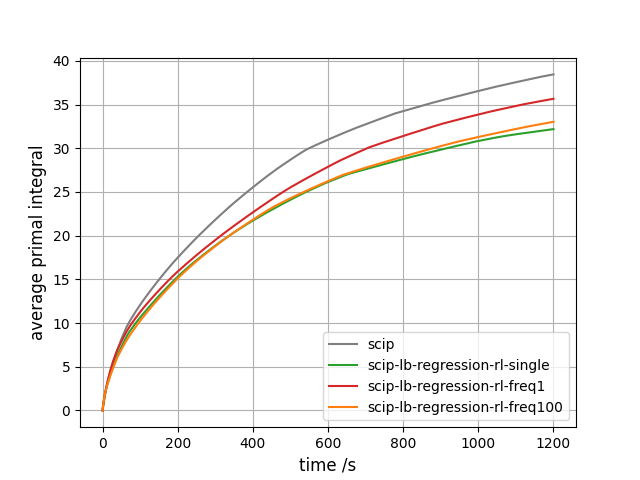}
    \caption{Evolution of \emph{primal integral} (geometric means) over time on binary MIPLIB dataset (1200s).}
    \label{fig:pi_new_long_multicall}
\end{figure}

We evaluate the two primal heuristic implementations on the MIPLIB binary dataset and report the evolution of the average \emph{primal integral} in Figure \ref{fig:pi_new_long_multicall}. Our observation is that the SCIP solver is improved by adding our learning based LB into SCIP as a primal heuristic. The results also suggest that the best strategy of executing LB is only to call it at the root node of the B\&B tree (or at the node where the first incumbent solution is found).


\section{Discussion}
\label{sec:discussion}
In this work, we have looked at the local branching paradigm by using a machine learning lens. We have considered the neighborhood size as a main factor for quantifying high-quality LB neighborhoods. We have presented a learning based framework for predicting and adapting the neighborhood size for the LB heuristic. The framework consists of a two-phase strategy. For the first phase, a \emph{scaled regression} model is trained to predict the size of the LB neighborhood at the first iteration through a regression task. In the second phase, we leverage reinforcement learning and devise a \emph{reinforced neighborhood search} strategy to dynamically adapt the size at the subsequent iterations. Furthermore, we have also investigated the relation between ``LB node" time limit $t$ and the actual performance of LB scheme, and devised a strategy for adapting $t$.
We have computationally shown that the neighborhood size and LB node time limit can indeed be learned, leading to improved performances and that the overall algorithm generalizes well both with respect to the instance size and, remarkably, across instances. 

Our framework relies on the availability of an initial solution, thus it can be integrated with other refinement heuristics. For future research, it would be  interesting to design more sophisticated hybrid frameworks that learn to optimize multiple refinement heuristics in a more collaborative way.



\begin{acknowledgements}
We would like to thank Maxime Gasse, Didier Chételat and Elias Khalil for their helpful discussions on the project. This work was supported by Canada Excellence Research Chair in Data Science for Real-Time Decision-Making at Polytechnique Montréal. The work of Matteo Fischetti was partially supported by MiUR, Italy.
\end{acknowledgements}

%
%


\bibliographystyle{spmpsci}      

\bibliography{main}

%
%



\section*{Statements and Declarations}

\paragraph{\bf Funding} This work was supported by Canada Excellence Research Chair in Data Science for Real-Time Decision-Making at Polytechnique Montréal. The work of Matteo Fischetti was partially supported by MiUR, Italy.

\paragraph{\bf Conflict of interest} The authors declare that have no conflict of interest.

\paragraph{\bf Competing Interests} The authors have no relevant financial or non-financial interests to disclose.

\paragraph{\bf Author Contributions}
All authors contributed to the study and development of this work. All authors read and approved the final manuscript.

\paragraph{\bf Data Availability} The datasets used the current study are either public, or available from the corresponding author on reasonable request.

\paragraph{\bf Code Availability} The code is publicly available at \url{https://github.com/pandat8/ML4LB}

\end{document}